\documentclass{article}
\usepackage[russian,english]{babel}
\usepackage{amssymb,amsmath,amsfonts}
\usepackage{graphics,dcolumn,bm,epic,eepic,float}
\usepackage{graphicx} 
\usepackage{epstopdf} 
\usepackage{cite}

\usepackage{bm} 
\usepackage{subfigure}
\usepackage{color} 

\usepackage{algorithm2e}
\usepackage[pdftex]{hyperref}

\usepackage{color, colortbl}
\definecolor{Yellow}{rgb}{1,1,0}
\definecolor{Grey}{rgb}{.87,.87,.87}
\definecolor{Purple}{rgb}{.8,.0,1.0}

\usepackage[english]{babel}
\newtheorem{theorem}{\hspace{\parindent}Theorem}
\newtheorem{lemma}{\hspace{\parindent}Lemma}

\begin{document}
\title{Existence of Nontrivial Negative Resonances for Polynomial 
Ordinary Differential Equations With Painlev\'e Property}

\date{}

\author{Stanislav Sobolevsky
	\thanks{SENSEable City Laboratory, Massachusetts Institute of Technology, 77 Massachusetts Avenue 9-221, Cambridge, MA 02139, USA}
	\thanks{E-mail: stanly@mit.edu}
}

\maketitle

\begin{abstract}
The Painlev\'e classification is one of the central problems in analytics theory of differential equations rooted in the XIX century. Although it saw many significant advances in analyzing certain classes of equations, the classification still remains an open problem especially for the higher-order equations. One of the main classical methods of Painlev\'e analysis is based on considering the resonance numbers corresponding to the possible indices of arbitrary coefficients in the Laurent expansion of the general solution in a neighborhood of a movable singularity. Complex and non-integer values of resonance numbers point out to existence of the movable critical singularities and positive integer numbers could be used to construct the said general solution. Also the equation always possesses at least one negative resonance number of $-1$ which corresponds to an arbitrary position of a movable pole. However our understanding of the role of nontrivial negative resonances different from $-1$ remains limited in spite of certain recent methodological advances related to it. And though in the lower-order classifications built so far such equations with nontrivial negative resonances have rather been a special case, the result of present work demonstrates that negative resonances are in fact common for the higher degree ordinary differential equations with Painlev\'e property. Specifically we'll prove that their presence is the necessary condition of the Painlev\'e property for the equations with degree of the leading terms higher than two.
\end{abstract}

\section*{Introduction}
The problem of Painlev\'e classification of the ordinary differential equations (i.e. finding all the classes of equations free of movable critical singularities) saw significant advances over the last century since it has been formulated. For the differential equations of the first order algebraic in independent variable and it's derivative the classification is based on the well-known Fuchs conditions (for the review see for instance \cite{ince1956}, chapter XIII). For the
second order equations with rational right-hand side the problem has been solved in the classical
works of Painlev\'e and Gambier (for the review see for instance \cite{ince1956},
chapter XIV). For the third order equations with polynomial right-hand side the classification has been started in the famous work of Chazy \cite{chazy1911equations} and recently
completed by C.Cosgrove \cite{cosgrove2000chazy}. Finally for the fourth order the problem has been solved by C.Cosgrove \cite{cosgrove2000P2, cosgrove2006P1} - again just in the polynomial case. The higher-order problem still remains unsolved in the more or less general case although certain special cases of the higher-order equations have been successfully classified in recent works, including binomial equations \cite{sobolevsky2005binomial3, sobolevsky2006binomialN}, nonlinear algebraic differential equations of arbitrary order $n$ missing two derivatives of order $n-1$ and $n-2$ \cite{sobolevskii2003algsing,   sobolevsky2004stam, sobolevskii2005alg, sobolevsky2006mono} and second-degree equations of arbitrary order \cite{sobolevsky2014quad}. 

One of the central classical methods of Painlev\'e analysis often called Painlev\'e test is based on studying the Laurent series representing the general solution around the possible movable singularity. The indices of the possible arbitrary coefficients called resonances are used as an important characteristic of the equation, while having all integer values for the resonances is used as a necessary condition for the Painlev\'e property. One of the resonances is always $-1$ corresponding to an additional arbitrary parameter being the location of the considered singularity. The other integer values of resonances if positive are used in order to construct the Laurent series for the solution. Requirement of it's consistency further leads to additional necessary conditions. However negative values of resonances are not involved in the classical version of the Painlev\'e test \cite{ablowitz1980connection}. Further some partial methodologies were presented to extract additional valuable information about the equation from negative resonances, including detecting logarithmic branch points \cite{fordy1991analysing, conte1993perturbative, musette1995non}. However so far our understanding of the role of negative resonances remains limited and a single consistent methodology of their utilization is still missing. Also equations with nontrivial (other than $-1$) negative resonances so far appeared within the classifications of the lower-order differential equations rather as a special case. In the present work we show that this case of negative resonance is actually common for the higher-order higher-degree equations having Painlev\'e property.

In \cite{Sobolevsky2004poly, sobolevskii2006modification} it has been proved that each nonlinear polynomial differential equation with the Painlev\'e property possesses movable singularities and, in particular, movable poles. In \cite{sobolevsky2014quad} it has been demonstrated that if the equation has a strong Panlev\'e property (all movable singularities are poles) then at least one of such movable poles should possess a complete set of $n-1$ positive resonances with the only trivial negative resonance number of $-1$. If the degree of the right-hard side is just $2$ then this pole is the only possible singularity of the equation and all such equations were classified \cite{sobolevsky2014quad}. In this paper we consider the polynomial ordinary differential equations of higher degree (i.e. having the degree of the leading terms higher than $2$), proving that such equations with Painlev\'e property always contain other movable poles having nontrivial negative resonance numbers as well.

\section{Main definitions}

Consider the non-linear polynomial ordinary differential equation
\begin{equation}
w^{(n)}=P(w^{(n-1)},w^{(n-2)},...,w,z),\label{pol}
\end{equation}
where $P$ is a polynomial in $w$ and its derivatives with coefficients locally
analytic in $z$. For the considered equations (\ref{pol}) similarly to \cite{sobolevsky2014quad} present the definition 
of the Painlev\'e property and the strong Painleve property. They correspond to the commonly accepted (but also often mixed with each other) concepts of freedom from movable branch points and from movable non-polar singularities correspondingly. 

\textbf{Definition\,1.} \cite{sobolevsky2014quad}
\begin{it} Consider an arbitrary solution
$w=w(z)$ of the equation (\ref{pol}) analytic in the
neighborhood of some point $z^*\in U$ and a path $\Gamma$ with the
beginning in $z^*$ along which all the coefficients of (\ref{pol})
can be analytically continued while the analytical continuation of
$w(z)$ comes to a singularity. Such a singularity of $w(z)$ is
called movable singularity of the considered solution $w=w(z)$.
\end{it}

\textbf{Definition\,2.} \cite{sobolevsky2014quad}
\begin{it}
The equation (\ref{pol}) is called to possess the Painlev\'e
property if the equation's solutions are single-valued near all of their movable
singularities. The equation (\ref{pol}) is called to possess
the strong Painlev\'e property when all (if any) movable singularities of it's
solutions are poles.
\end{it}

The definitions above are pretty much aligned with the classical understanding \cite{ince1956,conte1999painleve} of the concept just 
specifying it's technical details in the mathematically rigorous form. While definition of the Painlev\'e property corresponds to the most commonly accepted notion of freedom from movable branch points, the equations with the strong Painlev\'e property according to the definition 2 have also been introduced in \cite{gordoa2003new} as the equations of Painlv\'e-type.

\section{General conditions of Painlev\'e property}

In this section similarly to \cite{Sobolevsky2004poly, sobolevsky2014quad} we present a formalism for the well-known technique of the Painlev\'e test \cite{ablowitz1980connection} determining if the considered equation could possess the Painlev\'e property. 

If the equation (\ref{pol}) admits a solution with movable
pole of a certain order $s$ in a certain point $z=z_0$, there should exist 
a Laurent expansion of the form
\begin{equation}
w=q(z-z_0)^{-s}+\sum\limits_{j=1}^\infty c_j (z-z_0)^{j-s}, \label{laurent}
\end{equation}
converging in a certain deleted neighborhood of $z=z_0$,
where $s$ is an integer order of the pole, while $q, c_j$ are complex coefficients, provided that $q\neq 0$. 



Substituting series (\ref{laurent}) with undefined coefficients
$q, c_j$ into the equation (\ref{pol}) or (\ref{eqmain}) one should obtain an
algebraic equation
\begin{equation}
\label{determining}
H(q,z_0)=0
\end{equation}
for defining the possible values of major coefficient $q$ (so called determining equation) and
the following recurrent equations for determining each further coefficient $c_j$ through
the earlier defined $q,c_1,...,c_{j-1}$:
\begin{equation}
R(j,q,z_0)q_j=Q_j(q,c_1,\ldots,c_{j-1},z_0), \label{ajform}
\end{equation}
where $H$, $R$ and $Q_j$ are polynomials in $q$ and $c_j$. 


One should notice however that the coefficients of equations (\ref{determining}) and (\ref{ajform}) depend on $z_0$. 
For any $z_0$ different from the zeros of the coefficients of the original equation (\ref{pol}) the determining equation (\ref{determining}) and the equations (\ref{ajform}) have the same 
general form, while coefficients of the polynomials $H,R,Q_j$ are analytic in $z_0$ being polynomially expressed through the values of original equation's coefficients and their derivatives in $z_0$\cite{sobolevsky2014quad}.


For each root $q$ of the determining equation (\ref{determining}), one could use the equations (\ref{ajform})
to find all of the other coefficients $c_j$ of (\ref{laurent}) except of those who's indices $j=r$ are
the roots of the equation
\begin{equation}\label{resonance}
R(r,q,z_0)=0,
\end{equation}
called resonance equation, corresponding to the given choice of $q$. Its roots
are called resonance numbers (or simply resonances) of the equation (\ref{pol}) for the
selected $q$. 

Now going back to the considered equation (\ref{pol}) rewrite it in the following form
\begin{equation}
\begin{array}{c}
w^{(n)}=\sum\limits_{\chi\in \Omega}a_\chi(z)\prod\limits_{j=0}^{n-1}\left(w^{(j)}\right)^{\chi_j},
\end{array}
\label{eqmain}
\end{equation}
where $\chi$ are multi-indexes $\chi=(\chi_0, \chi_1,\ldots , \chi_{n-1})$ from a certain set $\Omega$ and $a_\chi$ are the equation's coefficients (not identically zero) analytic in $z$ in a certain complex domain.
Introduce a constant ${\cal B}=\min\limits_{\chi\in \Omega, |\chi|>1} \frac{n-\nu(\chi)}{|\chi|-1}$, where $|\chi|=\sum\limits_{j=0}^{n-1}\chi_j$ denotes the term's degree and $\nu(\chi)=\sum\limits_{j=0}^{n-1}j\chi_j$ --- the total number of derivative operators it contains. This characteristic is called the Bureau number \cite{Sobolevsky2004poly} and indicates the typical order $s$ of a movable pole that the solutions of equation (\ref{eqmain}) could contain. Although sometimes equation (\ref{eqmain}) could admit certain partial solutions having movable poles of a different order, like does a well-known example of Bureau equation \cite{Bureau1964}
$$
w^{(IV)}+3ww''-4\left(w'\right)^2=0
$$
having ${\cal B}=2$, but admitting a family of partial solutions $w=C(z-z_0)^{-3}-60(z-z_0)^{-2}$ with arbitrary $z_0$ and $C$ and a movable pole of order $3$. This equation however does not possess a Painlev\'e property.

According to the theorem 4 \cite{Sobolevsky2004poly} if the initial equation (\ref{eqmain}) possesses
the Painlev\'e property, the number ${\cal B}$ could only be $1$ or $2$. Moreover with respect to the theorem 3 \cite{Sobolevsky2004poly} at least one of the leading coefficients $a_{(1,0,0,....,0,1)}(z)$, $a_{(0,1,0,....,0,1,0)}(z)$, $a_{(2,0,0,....,0,1,0)}(z)$ in case ${\cal B}=1$ or the coefficient $a_{(1,0,0,....,0,1,0)}(z)$ in case ${\cal B}=2$ is not identically equal to zero as the equation's learning terms should include $w^{(n-1)}$ or $w^{(n-2)}$.

Now for the equation in the form (\ref{eqmain}), looking for a movable pole of a typical order $\cal B$ in a certain location $z_0$, one can specify the determining (\ref{determining}) and resonance (\ref{resonance}) equations in the following forms:
\begin{equation}
\begin{array}{c}
0=H(q,z_0)=\tau(-s,n)q-\sum\limits_{\chi\in \Omega_0}a_\chi(z_0)\left(\prod\limits_{j=0}^{n-1}\tau(-s,j)^{\chi_j}\right)q^{|\chi|},
\end{array}
\label{detEq}
\end{equation}
\begin{equation}
\begin{array}{c}
R(r,q,z_0)=\tau(-s+r,n)-\sum\limits_{\chi\in \Omega_0}a_\chi(z_0)\left(\prod\limits_{j=0}^{n-1}\tau(-s,j)^{\chi_j}\right)q^{|\chi|-1}\sum\limits_{j=0}^{n-1}\chi_j\frac{\tau(-s+r,j)}{\tau(-s,j)}
\end{array}
\label{resEq}
\end{equation}
where $\tau(-s,j)=\prod\limits_{k=0}^j(-s-k)$ and $\Omega_0=\{\chi\in \Omega: |\chi|+{\cal B}\nu(\chi)=n+{\cal B}\}$ is the set of equation's leading terms. Let $d=\max\limits_{\chi\in \Omega_0}|\chi|$ be the maximal degree of those terms and notice that it as well as $\Omega_0$ does not depend on the choice of $z_0$.

Also let $m={\rm deg}_q H(q,z_0)-1$ be the order of the determining equation (\ref{detEq}) after it's reduction through division by $q$. Now denote the $m$ non-zero roots (including multiple ones) of the corresponding equation (\ref{detEq}) by $q^1,q^2,\ldots$, $q^m$ (as one can see a constant term of the polynomial $H/q$ is nonzero, i.e. the $m$ roots are nonzero too). Let $\hat{h}\neq 0$ be the major coefficient of the polynomial $H(q,z_0)/q$. Then for each $q=q^k$ consider the set of the corresponding roots $r^k=(r^k_1,r^k_2,\ldots,r^k_n)$ of the resonance equation (\ref{resEq}). A well known necessary condition for the Painlev\'e property (in particular rigorously proved in \cite{Sobolevsky2004poly}) claims that for each $k$ all the resonance numbers $r^k_1,r^k_2,\ldots,r^k_n$ should be distinct integers.

The easiest way of proving the condition above is by introducing into the equation (\ref{eqmain}) a small parameter transform $z=z_0+\alpha Z$, $w=\alpha^{-{\cal B}} W$ which for the limit case of $\alpha=0$ gives a reduced equation
\begin{equation}
\begin{array}{c}
W^{(n)}=\sum\limits_{\chi\in \Omega_0}a_\chi(z_0)\prod\limits_{j=0}^{n-1}\left(W^{(j)}\right)^{\chi_j},
\end{array}
\label{eqreduced}
\end{equation}
with parametric families of partial solutions
\begin{equation}
W=W_0(Z)=q (Z-Z_0)^{-{\cal B}},
\label{solW0}
\end{equation}
were $q=q^k$ is an arbitrary root of (\ref{detEq}). Obviously determining (\ref{detEq}) and resonance (\ref{resEq}) equations for the reduced equation (\ref{eqreduced}) are exactly the same as for the initial equation (\ref{eqmain}) for the considered $z_0$.
Then solutions of the transformed initial equation for $\alpha\neq 0$ on a certain circle $\Gamma$ surrounding a singularity $Z=Z_0$ could be expressed by means of a converging series $W(Z)=W_0(Z)+\sum\limits_{j=1}^\infty \alpha^j W_j(Z)$, where all functions $W_j$ are analytic on $\Gamma$. Also in this case $W_1$ satisfies a homogenous Euler linear differential equations with a general solution $W_1=\sum_{i=1}^n C_i (Z-Z_0)^{-{\cal B}+r^k_i}$ if all roots $r^k_i$ are distinct, taking a more complex form with a logarithmic singularity in $Z=Z_0$ otherwise.

However the resonance numbers $r^k_1,r^k_2,\ldots,r^k_n$ being integer could be of both signs - positive and negative. Furthermore for any $k$ one of the resonance numbers is always equal to $-1$ (see \cite{ablowitz1980connection} or more recently \cite{conte1999painleve} (page 126-127)). In fact this could be easily proven by looking at the partial derivative of the equation (\ref{eqreduced}) by $Z_0$ along the above solution (\ref{solW0}); this shows that $\partial W_0(Z)/\partial Z_0=-{\cal B}q (Z-Z_0)^{-{\cal B}-1}$ should satisfy the above Euler linear differential equation.

Finally notice that none of the resonance numbers $r^k_1,r^k_2,\ldots,r^k_n$ could be equal to zero as otherwise $W=C (Z-Z_0)^{-{\cal B}}$ could be proved to satisfy (\ref{eqreduced}) for an arbitrary $C$ which could not be the case as only certain values of the coefficient satisfy the determining equation (\ref{detEq}) for $m>0$.


According to Vieta's theorem for the resonance equation (\ref{resEq}) the product $Pr(q^k,z_0)=\prod\limits_{i=1}^n r^k_i$ of resonance numbers $r^k$ corresponding to each $k$ could be found as 
$$
Pr(q^k,z_0)=\prod\limits_{j=1}^{n} r_j=(-1)^n R(0,q^k,z_0)=(-1)^n\frac{\partial H}{\partial q}(q^k,z_0)=
$$$$
=\hat{h}q^k(-1)^n\prod\limits_{j, j\neq k}(q^k-q^j).
$$
In turn according to Vieta's theorem for the reduced determining equation $H(q,z_0)/q=0$ one has $\prod\limits_j q^j=(-1)^{n+m}(n+{\cal B}-1)!/\hat{h}$ and this way $\hat{h}=(-1)^{m+n}(n+{\cal B}-1)!/\prod\limits_j q^j$. Consequently,
\begin{equation}
Pr(q^k,z_0)=-(n+{\cal B})!\prod\limits_{j, j\neq k}\frac{q^j-q^k}{q^j}.
\label{resProds}
\end{equation}
In particular it means that the determining equation (\ref{detEq}) is free from multiple roots as otherwise if $q=q^k$ is the one, then $Pr(q^k,z_0)=0$ which is not possible as the resonances are nonzero numbers.

Further as one can see the conditions (\ref{resProds}) lead to the following equation
\begin{equation}
\sum\limits_{k=1}^m \frac{1}{Pr(q^k,z_0)}=-\frac{1}{(n+{\cal B}-1)!}.
\label{resProds2}
\end{equation}
Indeed, with respect to (\ref{resProds}) the above equation is equivalent to
$$
\sum\limits_{j=1}^m \frac{1}{\partial H(q^j,z_0)/\partial q}=-\frac{1}{\partial H(0,z_0)/\partial q}.
$$ 
Since all the roots of $H(q,z_0)$ are non-multiple, one can get the last equation following from the residual theorem for the function $1/H(q,z_0)$, considered as the function of $q$ in complex domain.
The condition (\ref{resProds2}) was obtained for the lower-order equation in the previous works \cite{chazy1911equations, Bureau1964, martynov1973differential, cosgrove2006P1}.

Another condition which could be also easily obtained from the Vieta's theorem for (\ref{resEq}) is related to the sum of the resonance numbers \cite{Sobolevsky2004poly}. If ${\cal B}=1$ one can get for any $k$
\begin{equation}
\sum\limits_{j=1}^n r^k_j=\sum\limits_{j=1}^n j+A q^k,
\label{resSums}
\end{equation}
where $A=a_{(1,0,0,\ldots,0,1)}(z_0)$, while if ${\cal B}=2$,
\begin{equation}
\left\{\begin{array}{l}
\sum\limits_{j=1}^n (r^k_j)=\sum\limits_{j=1}^n j,\\
\sum\limits_{j=1}^n (r^k_j)^2=\sum\limits_{j=1}^n j^2+B q^k,
\end{array}\right.
\label{resSums2}
\end{equation}
where $B=a_{(1,0,\ldots,0,1,0)}(z_0)$. Also notice that in the last case for the equation (\ref{eqmain}) with Painlev\'e property and ${\cal B}=2$, from the first equation of (\ref{resSums2}) one can see that $\sum\limits_{j=1}^n (r^k_j)^2>\sum\limits_{j=1}^n j^2$ given that the set $r^k$ consists of distinct integers (proof of this statement could be found in \cite{sobolevskii2005alg}). So from the last equation of (\ref{resSums2}) all $B q^k$ should be positive real numbers.

The conditions above present some considerable restrictions on the important characteristics of the equation - values of $q^1,q^2,...,q^m$ and sets of $r^1,r^2,...,r^m$. However the structure of these characteristics first of all depends on the value of $m$. And although based on the structure of (\ref{detEq}), $m$ obviously appears to be related to the maximal degree of the equation's leading terms $d$ ($m\leq d-1$), but in some special cases when some of the terms of (\ref{detEq}) vanish, $m$ could be actually lower than $d-1$. For example for the equation 
$$
w'''=w''w-2(w')^2
$$
having Bureau number $1$ and degree $d=2$, all terms of (\ref{detEq}) but $-6q$ vanish and as a result we have $-6q=0$ and $m=0$. 

However according to \cite{sobolevskii2006modification} equations of this special kind do not possess the Painlev\'e property. We have to notice that in \cite{sobolevskii2006modification} this is formally claimed for the nonlinear reduced equation (\ref{eqreduced}) in case $m=0$ only, i.e. once all the terms of equation (\ref{detEq}) vanish except $\tau(-s,n)q$, while the equation (\ref{eqreduced}) remains nonlinear. However the proof of this claim (theorem 2 \cite{sobolevskii2006modification}) in fact only relies on the condition 
\begin{equation}
\sum\limits_{\chi\in \Omega_0, |\chi|=d}a_\chi(z_0)\left(\prod\limits_{j=0}^{n-1}\tau(-s,j)^{\chi_j}\right)=0,
\label{vanish}
\end{equation}
which is the one sufficient for the condition (29) \cite{sobolevskii2006modification} to hold. This way each time $m<d-1$ which means that the highest degree terms vanish, i.e. the condition (\ref{vanish}) holds, the nonlinear equation (\ref{eqreduced}) as well as the initial equation (\ref{eqmain}) does not possess the Painlev\'e property. Of course this does help to deal with cases when $z_0$ is the common root of the coefficients of all the leading terms but $w^{(n)}$ (having a constant coefficient) - in such a case (\ref{eqreduced}) is a trivial linear equation. But other than that if $z_0$ does not belong to the countable set $E$ of isolated points being the roots of equation's coefficients, one finds $m$ to be equal to $d-1$ for the equation (\ref{eqmain}) to have the Painlev\'e property. Further consider $z_0\not\in E$.

\section{Existence of negative resonances}
Now opposed to the statement we want to prove, assume that for any $j$ the set of resonance numbers $r^j$ does not contain any non-trivial negative values other than $-1$. Then $Pr(q^j,z_0)<0$ for all $j$. Consider several possible cases for the roots $q^j$ of determining equation (\ref{determining}), without loss of generality assuming that $q_1=1$ (otherwise one could always achieve it introducing a change of variables $w=q_1 W$):\\
A) All roots of (\ref{determining}) are real numbers, while $m>2$.\\
B) Set of roots (\ref{determining}) contain complex numbers, while $m>2$.\\
C) Number of roots $m=2$.

In case A) there always exist two different roots $q^i<q^j$ of the same sign (i.e. either both positive, either both negative real numbers) and without loss of generality one can consider that an interval $(q^i,q^j)$ does not contain any other roots $q^k$. Then according to (\ref{resProds2}) as one can see $Pr(q^i,z_0)Pr(q^j,z_0)<0$ (as changing $i$ to $j$ preserves the sign of all terms of  (\ref{resProds}) but $(q^i-q^j)$) which contradicts an assumption above.

The case B) is impossible for ${\cal B}=2$ because of (\ref{resSums2}) and $B\neq 0$. While if ${\cal B}=1$ then due to (\ref{resSums}) the only possible option is $A=0$. However in this case for any $i$, we have $\sum\limits_{k=1}^n r^i_k=\sum\limits_{k=1}^n k$. Then as all resonance numbers but $-1$ are positive integers, no single set $r^i$ can contain $1$ as otherwise with respect to (\ref{resSums}) one can see that $|Pr(q^i,z_0)|< n!$. Indeed if we denote $r^i_1=-1, r^i_2=1$, then $\sum\limits_{k=3}^n r^i_k=n(n+1)/2$ and $|Pr(q^i,z_0)|=\prod\limits_{k=3}^n r^i_k\leq n!$ as a consequence of the following

\begin{lemma}
\label{lemP}
For the product $P$ of any set of $t$ pairwise different natural numbers $y_1,y_2,...,y_t$ the following inequality holds:
\begin{equation}
P\leq P^{max}(t,S)=\frac{(t+\tau)!}{(\tau-1)!\zeta},
\label{maxProd}
\end{equation}
where $S=\sum\limits_{i=1}^t y_i$ is the sum of the given numbers, $\zeta=\tau+(1-\epsilon)t$, while $\tau$ and $\epsilon$ are integer and fractional parts of $S/t-(t-1)/2$, i.e. $\tau=\left[S/t-(t-1)/2\right]$, $\epsilon=\left\{S/t-(t-1)/2\right\}$.
\end{lemma}

In plain language the lemma above claims that the highest possible product of a set of the given amount of distinct natural numbers with the given sum is achieved for the most dense distribution of said numbers, i.e. $\{\tau,\tau+1,...,\tau+t\}\setminus\{\zeta\}$ (while $\tau$ and $\zeta$ are defined in a way consistent with the given $S$). The lemma could be easily proved using the method of mathematical induction. I.e. being obvious for $t=1$, if validated for all sets of a certain size $t=k$ this claim then follows for the sets of the higher size of $t=k+1$. Indeed assume that $y_1<y_2<...<y_t$ and let $y_1=l<\tau$ (as if $y_1\geq \tau$ then the set $\{\tau,\tau+1,...,t+\tau\}\setminus\{\zeta\}$ is the only possible option for $y_1,y_2,...,y_t$ and so $P=P^{max}(t,S)$). Then using the claim of the lemma for $y_2,y_3,...,y_t$ with the sum $S-l$ one can show that $l P^{max}(t-1,S-l)<(l+1)P^{max}(t-1,S-l-1)<...<(\tau-1)P^{max}(t-1,S-\tau+1)<\tau P^{max}(t-1,S-\tau)=P^{max}(t,S)$ since $P^{max}(t-1,S-(x-1))/P^{max}(t-1,S-x)=\zeta^*/(\zeta^*-1)<\tau/(\tau-1)\leq x/(x-1)$ for any $x\leq \tau$, where 
$\zeta^*=\tau^*+(1-\epsilon^*)t>\tau$, while $\tau^*=\left[(S-x)/(t-1)-(t-2)/2\right]>\tau$, $\epsilon^*=\left\{(S-x)/(t-1)-(t-2)/2\right\}$.

Going back to the set $r^i$ as described above we get $\tau=3$ for it in case $n>5$ or $\tau=4$ otherwise. Then according to the lemma \ref{maxProd}, for $n>5$ we have $P\leq (n+1)!/(n-2)/2<n!$, for $n=5$ and $P\leq 6!/3!=n!$, for $n=4$, $P\leq 6!/3!/5=n!$, which completes the proof of the claim above.

Then as $\forall j, Pr(q_j,z_0)<0$ we have 
$$
\sum\limits_{j=1}^m \frac{1}{Pr(q^j,z_0)}<\frac{1}{Pr(q^i,z_0)}\leq -\frac{1}{n!}
$$
which comes into a contradiction with (\ref{resProds2}).

Now if every set $r^i$ consists of $-1$ and integers greater than $1$, then for all such sets due to (\ref{resSums}) only two options are possible: $r^i=(-1,2,3,\ldots,n-1,n+2)$ or $r^i=(-1,2,3,\ldots,n-2,n,n+1)$. In both cases for $n\geq 3$, according to the lemma \ref{maxProd} one can see that $|Pr(q_i,z_0)|\leq (n+1)!/(n-1)\leq 2 n!$ for any $i$, which once again comes into contradiction with (\ref{resProds2}).

This way the only possible case is C), i.e. $m=2$. In this case we have $Pr(q_1,z_0)=(q_2-1)(n+{\cal B}-1)!/q_2$ and $Pr(q_2,z_0)=(1-q_2)(n+{\cal B}-1)!$. This way $q_2$ is a real number. If $q_2>0$ then at least one of $Pr(q_1,z_0)$, $Pr(q_2,z_0)$ should be positive, which comes into a contradiction with assumption above. Consequently the only possible option is $q_2<0$, which due to (\ref{resSums2}) is possible only for ${\cal B}=1$.

Then $|Pr(1,z_0)|>n!$ and $|Pr(q_2,z_0)|>n!$ and because of (\ref{resSums}) for at least one of $i=1,2$, the inequality $\sum\limits_{k=1}^n r^i_k\leq \sum\limits_{k=1}^n k$ holds. Without loss of generality assume that this is $i=1$ as otherwise it is sufficient to introduce a variable transform $w=q_2 v$.

Then $r^1$ does not contain $1$ as otherwise $|Pr(1,z_0)|<n!$. 
However this is possible only in one of the following four cases:\\
1) $r^1=(-1,2,3,\ldots,n)$;\\
2) $r^1=(-1,2,3,\ldots,n-1,n+1)$;\\
3) $r^1=(-1,2,3,\ldots,n-2,n,n+1)$;\\
4) $r^1=(-1,2,3,\ldots,n-1,n+2)$;\\

First case is not possible as this way $|Pr(1,z_0)|=n!$. 

In the second case $|Pr(1,z_0)|=n!\frac{n+1}{n}=n!\frac{q_2-1}{q_2}$, i.e. $q_2=-n$, $|Pr(q_2,z_0)|=(n+1)!$,
and $\sum r^2_j =\sum\limits_{j=1}^n j+n$.
If $r^2$ contains $r=1$ then according to lemma \ref{lemP}, $|Pr(q_2,z_0)|\leq (n+2)!/(6(n-2))<(n+1)!$ 
for any $n\geq 3$.
If $r=2$ is contained in $r^2$ then according to lemma \ref{lemP}, $|Pr(q_2,z_0)|\leq (n+2)!/(3(n-1))<(n+1)!$ for any $n\geq 3$.
Otherwise if $r=1,2$ are not contained in $r^2$ then because of $\sum\limits_{j=1}^n r^2_j =\sum\limits_{j=1}^n j+n$ the only possible options are:\\
2a) $r^2=(-1,3,4,\ldots,n,n+4)$;\\
2b) $r^2=(-1,3,4,\ldots,n-1,n+1,n+3)$;\\
2c) $r^2=(-1,3,4,\ldots,n-2,n,n+1,n+2)$.\\
In case 2a) one has $(n+1)!=|Pr(q_2,z_0)|=n!\cdot (n+4)/2$, i.e. $n=2$.
In case 2b) one has $(n+1)!=|Pr(q_2,z_0)|=(n+1)!\cdot (n+3)/(2n)  $, i.e. $n=3$.
In case 2c) one has $(n+1)!=|Pr(q_2,z_0)|=(n+1)!\cdot (n+2)/(2(n-1))  $, i.e. $n=4$. However if $n=4$ then $r^2$ as defined contains $r=2$ which contradicts the assumption above. 

In cases 3) and 4) the condition $\sum\limits_{j=1}^n r^1_j= \sum\limits_{j=1}^n j$ holds, therefore with respect to (\ref{resSums}) one has $A=0$. Consequently $r^2$ in the same way is also one of the two sets $(-1,2,3,\ldots,n-2,n,n+1)$ or $(-1,2,3,\ldots,n-1,n+2)$. If sets $r^1$ and $r^2$ are different then with respect to (\ref{resProds2})
$$
\frac{n-1}{(n+1)n!}+\frac{n}{(n+2)n!}=\frac{1}{n!},
$$
but this does not hold for integer $n$. While if the sets $r^1$ and $r^2$ are the same, then in the case 3) one has
$$
\frac{2(n-1)}{(n+1)n!}=\frac{1}{n!},
$$
i.e. $n=3$, and in the case 4) one has
$$
\frac{2n}{(n+2)n!}=\frac{1}{n!},
$$
i.e. $n=2$. This way for $n>3$ the case C) is also impossible.

Finally, we proved the following

\begin{theorem} 
\label{theoremNegRes}
If the equation (\ref{pol}) of order $n>3$ possesses the Painlev\'e property and the highest degree $d$ of it's leading terms $\Omega_0$ is more than $2$, then for any location $z_0\not\in E$ (i.e. different from the roots of equation's coefficients) the determining equation (\ref{detEq}) has $m=d-1>1$ nonzero roots, and for at least one of them the corresponding set of resonance numbers contains non-trivial negative values other than $-1$.
\end{theorem}

%
%


\section{Conclusion}

In this paper we demonstrated that negative values of resonance numbers are quite a common case for the nonlinear differential equations with polynomial right-hand side (\ref{pol}) having Painlev\'e property - the only type of equations without non-trivial negative resonances are those with the maximal degree of leading terms $d=2$. The last case however seems to be rather straightforward - all equations of the second degree having strong Painlev\'e property are already considered in detail in \cite{sobolevsky2014quad}, and only one class of such equations was found for the higher order case $n>6$ (and this single class appeared to be easily linearizable). All said it means that all new equations with Painlev\'e property one could expect to find in the class (\ref{pol}) will always possess nontrivial negative resonances emphasizing the paramount importance of their understanding and building appropriate methods for their analysis. 


\bibliography{Painleve}

\begin{thebibliography}{10}
\providecommand{\url}[1]{\texttt{#1}}
\providecommand{\urlprefix}{URL }
\expandafter\ifx\csname urlstyle\endcsname\relax
  \providecommand{\doi}[1]{doi:\discretionary{}{}{}#1}\else
  \providecommand{\doi}{doi:\discretionary{}{}{}\begingroup
  \urlstyle{rm}\Url}\fi
\providecommand{\bibAnnoteFile}[1]{%
  \IfFileExists{#1}{\begin{quotation}\noindent\textsc{Key:} #1\\
  \textsc{Annotation:}\ \input{#1}\end{quotation}}{}}
\providecommand{\bibAnnote}[2]{%
  \begin{quotation}\noindent\textsc{Key:} #1\\
  \textsc{Annotation:}\ #2\end{quotation}}
\providecommand{\eprint}[2][]{\url{#2}}

\bibitem{ince1956}
Ince E (1956) \selectlanguage{russian} Ordinary differential equations.
\newblock Dover, New York.
\bibAnnoteFile{ince1956}

\bibitem{chazy1911equations}
Chazy J (1911) Sur les {\'e}quations diff{\'e}rentielles du troisi{\`e}me ordre
  et d'ordre sup{\'e}rieur dont l'int{\'e}grale g{\'e}n{\'e}rale a ses points
  critiques fixes.
\newblock Acta Mathematica 34: 317--385.
\bibAnnoteFile{chazy1911equations}

\bibitem{cosgrove2000chazy}
Cosgrove CM (2000) Chazy classes ix--xi of third-order differential equations.
\newblock Studies in applied mathematics 104: 171--228.
\bibAnnoteFile{cosgrove2000chazy}

\bibitem{cosgrove2000P2}
Cosgrove CM (2000) Higher-order painlev{\'e} equations in the polynomial class
  i. bureau symbol p2.
\newblock Studies in applied mathematics 104: 1--65.
\bibAnnoteFile{cosgrove2000P2}

\bibitem{cosgrove2006P1}
Cosgrove CM (2006) Higher-order painlev{\'e} equations in the polynomial class
  ii: Bureau symbol p1.
\newblock Studies in Applied Mathematics 116: 321--413.
\bibAnnoteFile{cosgrove2006P1}

\bibitem{sobolevsky2005binomial3}
Sobolevsky S (2005) Binomial-type ordinary differential equations of the third
  order.
\newblock Studies in Applied Mathematics 114: 1--15.
\bibAnnoteFile{sobolevsky2005binomial3}

\bibitem{sobolevsky2006binomialN}
Sobolevsky S (2006) Painlev{\'e} classification of binomial type ordinary
  differential equations of the arbitrary order.
\newblock Studies in Applied Mathematics 117: 215--237.
\bibAnnoteFile{sobolevsky2006binomialN}

\bibitem{sobolevskii2003algsing}
Sobolevskii S (2003) Movable singular points of ordinary differential equations
  with algebraic singularities of the right-hand side.
\newblock Differential Equations 39: 381--386.
\bibAnnoteFile{sobolevskii2003algsing}

\bibitem{sobolevsky2004stam}
Sobolevsky S (2004) Movable singularities of a class of nonlinear ordinary
  differential equations of arbitrary order.
\newblock Studies in Applied Mathematics 112: 227--234.
\bibAnnoteFile{sobolevsky2004stam}

\bibitem{sobolevskii2005alg}
Sobolevskii S (2005) Movable singular points of algebraic ordinary differential
  equations.
\newblock Differential Equations 41: 1146--1154.
\bibAnnoteFile{sobolevskii2005alg}

\bibitem{sobolevsky2006mono}
Соболевский С (2006) Подвижные особые точки решений обыкновенных
  дифференциальных уравнений.
\newblock Минск: БГУ (in Russian), 118 pp.
\bibAnnoteFile{sobolevsky2006mono}

\bibitem{sobolevsky2014quad}
Sobolevsky S (2014) Painleve classification of polynomial ordinary differential
  equations of arbitrary order and second degree.
\newblock arXiv preprint arXiv:14102649 .
\bibAnnoteFile{sobolevsky2014quad}

\bibitem{ablowitz1980connection}
Ablowitz MJ, Ramani A, Segur H (1980) A connection between nonlinear evolution
  equations and ordinary differential equations of p-type. i.
\newblock Journal of Mathematical Physics 21: 715--721.
\bibAnnoteFile{ablowitz1980connection}

\bibitem{fordy1991analysing}
Fordy A, Pickering A (1991) Analysing negative resonances in the painlev{\'e}
  test.
\newblock Physics Letters A 160: 347--354.
\bibAnnoteFile{fordy1991analysing}

\bibitem{conte1993perturbative}
Conte R, Fordy AP, Pickering A (1993) A perturbative painlev{\'e} approach to
  nonlinear differential equations.
\newblock Physica D: Nonlinear Phenomena 69: 33--58.
\bibAnnoteFile{conte1993perturbative}

\bibitem{musette1995non}
Musette M, Conte R (1995) Non-fuchsian extension to the painlev{\'e} test.
\newblock Physics Letters A 206: 340--346.
\bibAnnoteFile{musette1995non}

\bibitem{Sobolevsky2004poly}
Sobolevskii S (2004) Movable singular points of polynomial ordinary
  differential equations.
\newblock Differential Equations 40: 807--814.
\bibAnnoteFile{Sobolevsky2004poly}

\bibitem{sobolevskii2006modification}
Sobolevskii S (2006) On a modification of the small parameter method.
\newblock Differential Equations 42: 218--228.
\bibAnnoteFile{sobolevskii2006modification}

\bibitem{conte1999painleve}
Conte R (1999) The Painlev{\'e} property: one century later, volume~1.
\newblock Springer New York.
\bibAnnoteFile{conte1999painleve}

\bibitem{gordoa2003new}
Gordoa PR, Joshi N, Pickering A (2003) A new technique in nonlinear singularity
  analysis.
\newblock Publications of the Research Institute for Mathematical Sciences 39:
  435--449.
\bibAnnoteFile{gordoa2003new}

\bibitem{Bureau1964}
Bureau F (1964) Differential equations with fixed critical points.
\newblock Annali di Matematica pura ed applicata 64: 229--364.
\bibAnnoteFile{Bureau1964}

\bibitem{martynov1973differential}
Martynov I (1973) Differential equations with stationary critical
  singularities.
\newblock Differential Equations 9: 1368--1376.
\bibAnnoteFile{martynov1973differential}

\end{thebibliography}

\end{document}